\documentclass{article}
\usepackage{amssymb}
\usepackage{amsmath}

\input{tcilatex}
\begin{document}

\begin{center}
{\LARGE Compound Lucas Magic Squares\medskip }

{\large Ronald P. Nordgren\footnote{%
email: nordgren@rice.edu}}

Brown School of Engineering, Rice University\vspace{-0.05in}{\LARGE \medskip 
}
\end{center}

\noindent \textbf{Abstract. }We review a general parameterization of an
order-$3$ magic square derived by Lucas and we compound it to produce a
parameterized order-9 magic square. Sequential compounding to order $3^{\ell
}$ $\left( \ell =2,3,\ldots \right) $ also is treated. Expressions are found
for the matrices in the Jordan canonical form and the singular value
decomposition of the compound Lucas magic square matrices. We develop a
procedure for determining if an order-$n$ magic square may be natural (with
elements $0,1,\ldots ,n^{2}-1).$ This enables determination of numerical
values for parameters in natural compound Lucas magic squares. Also, we find
commuting pairs of compound Lucas matrices and formulas for matrix powers of
order-3 and order-9 Lucas matrices. A parameterization due to Frierson is
related to Lucas' parameterization and our results specialize to it,
complementing previous results.\vspace{-0.1in}

\section{Introduction}

In 1894, Lucas \cite{LUCA} derived a general parameterization of an order-3
magic square matrix in terms of three parameters. We compound his order-3
matrix to produce a parameterized magic square matrix of order 9 having six
parameters. In addition, we sequentially compound his order-3 matrix to
produce parameterized magic squares of order $3^{\ell }$ $\left( \ell
=2,3,\ldots \right) $. Using our previous formulas for compound squares \cite%
{NORD2,NORD3}, we obtain expressions for the matrices in the Jordan
canonical form (JCF) and the singular value decomposition (SVD) of compound
Lucas matrices.

We develop a necessary condition for an order-$n$ magic square matrix to be
natural (with elements $0,1,\ldots ,n^{2}-1)$ based on its Frobenius norm
and singular values. Application of this condition and other considerations
proves that a Lucas square of order $3^{\ell }$ $\left( \ell =1,2,\ldots
\right) $ is natural when its parameters are $\pm 3^{k}$ $\left(
k=0,1,\ldots ,2\ell -1\right) $. We determine the number of fundamental
(phases excluded) natural Lucas magic squares that can be constructed by the
sequential compounding procedure. We find 48 fundamental numerical
realizations of natural Lucas magic squares of order 9. Also, we find
expressions for pairs of commuting Lucas matrices and we derive formulas for
matrix powers of order-3 and order-9 Lucas matrices.

Another parameterization of an order-3 magic square given by Frierson \cite%
{FIER} in 1907 is related to Lucas' parameterization by a linear
transformation of parameters. Frierson compounds his order-3 square to
produce a parameterized magic square matrix of order 9. After specifying one
of the parameters, he gives six fundamental numerical realizations of his
order-9 square and we find another six. Frierson's approach is reviewed by
Loly and Cameron \cite{LOLY3} who extend his compounding construction to
magic square matrices of order 27. They obtain formulas for singular values
of compound Frierson matrices of order $3^{\ell }$ $\left( \ell =2,3,\ldots
\right) $ based on extrapolation of calculations with Maple and Mathematica.
Their formulas for singular values are identical to ours for compound Lucas
matrices.

For a general overview of magic squares and a discussion of compounding
them, see Pickover \cite{PICK}. For previous work on compounding, see
Frierson \cite{FIER}, Bellew \cite{BELL}, Chan and Loly \cite{CHAN},
Eggermont \cite{EGGE}, Rogers, et. \negthinspace al \cite{ROGE}, Nordgren 
\cite{NORD2,NORD3}, and Loly and Cameron \cite{LOLY3}. For previous work on
the JCF and SVD of magic square matrices, see Mattingly \cite{MATT}, Loly,
et. \negthinspace al \cite{LOLY1}, Nordgren \cite{NORD1}, and Cameron, et.
\negthinspace al \cite{CAM}. Additional references are given in all of these
publications.\smallskip

\noindent \textbf{Definitions.} We recall \cite{PICK} that the rows, columns
and two main diagonals of an order-$n$ \emph{magic square }matrix $\mathbf{M}%
_{n}$ sum to the \emph{summation index} $\mu _{n},$ i.e.%
\begin{align}
\mathbf{M}_{n}\mathbf{E}_{n}& =\mathbf{E}_{n}\mathbf{M}_{n}=\mu _{n}\mathbf{E%
}_{n},  \notag \\
\limfunc{tr}\left[ \mathbf{M}_{n}\right] & =\limfunc{tr}\left[ \mathbf{R}_{n}%
\mathbf{M}_{n}\right] =\mu _{n},  \label{MnMagic}
\end{align}%
where $\mathbf{E}_{n}$ is the all-ones matrix and $\mathbf{R}_{n}$ is the
matrix with ones on the cross diagonal and all other elements zero. We
restrict our attention to $\mathbf{M}_{n}$ with non-negative integer
elements. When $\mathbf{M}_{n}$ is a \emph{natural} magic square with
elements $0,1,\ldots ,n^{2}-1,$ the summation index is given by%
\begin{equation}
\mu _{n}=\frac{n}{2}\left( n^{2}-1\right) \!.  \label{Mun}
\end{equation}%
This equation also can be satisfied by certain non-natural magic squares as
shown in Section 3.

In a \emph{regular} (associative) magic square, each pair of centrosymmetric
elements adds to $2\mu _{n}/n$ and in an odd-order square the center element
is half this value, i.e.%
\begin{equation}
\mathbf{M}_{n}+\mathbf{R}_{n}\mathbf{M}_{n}\mathbf{R}_{n}=\frac{2\mu _{n}}{n}%
\mathbf{E}_{n}.  \label{Reg}
\end{equation}%
We begin with the basic order-3 Lucas magic square which is compounded to
Lucas magic squares of higher order in Sections 4 and 5.

\section{Order-3 Lucas Square}

\noindent \textbf{Parameterizations}. Lucas \cite{LUCA} derived a general
parameterization of an order-3 magic square which can be written as%
\begin{equation}
\mathbf{L}_{3}\left( c,v,y\right) =\left[ 
\begin{array}{ccc}
c+v & c-v-y & c+y \\ 
c-v+y & c & c+v-y \\ 
c-y & c+v+y & c-v%
\end{array}%
\right] \!.  \label{Luc3}
\end{equation}%
It satisfies the magic square conditions (\ref{MnMagic}) with $\mu _{3}=3c$
and it is regular. When $\mathbf{L}_{3}$ is natural, by (\ref{Mun}) and (\ref%
{Luc3}), we have 
\begin{equation}
\mu _{3}=3c=12,\quad \therefore c=4.  \label{mu3}
\end{equation}%
In addition, in order that the maximum element of natural $\mathbf{L}_{3}$
equals $8$ and the minium element equals $0,$ the parameters $v$ and $y$
must be 
\begin{equation}
v=\pm 3,\ y=\pm 1\text{ or }v=\pm 1,\ y=\pm 3.  \label{vy3}
\end{equation}%
Thus, the following realizations of $\mathbf{L}_{3}\left( c,v,y\right) $ are
the only natural magic squares of order 3:%
\begin{gather}
\mathbf{L}_{3}\left( 4,3,1\right) =\left[ 
\begin{array}{ccc}
7 & 0 & 5 \\ 
2 & 4 & 6 \\ 
3 & 8 & 1%
\end{array}%
\right] \equiv \mathbf{\hat{L}}_{3},\quad \quad \quad \mathbf{L}_{3}\left(
4,3,-1\right) =\mathbf{\hat{L}}_{3}^{T}=\left[ 
\begin{array}{ccc}
7 & 2 & 3 \\ 
0 & 4 & 8 \\ 
5 & 6 & 1%
\end{array}%
\right] \!,  \notag \\
\mathbf{L}_{3}\left( 4,-3,-1\right) =\mathbf{R}_{3}\mathbf{\hat{L}}_{3}%
\mathbf{R}_{3}=\left[ 
\begin{array}{ccc}
1 & 8 & 3 \\ 
6 & 4 & 2 \\ 
5 & 0 & 7%
\end{array}%
\right] \!,\quad \mathbf{L}_{3}\left( 4,-3,1\right) =\mathbf{R}_{3}\mathbf{%
\hat{L}}_{3}^{T}\mathbf{R}_{3}=\left[ 
\begin{array}{ccc}
1 & 6 & 5 \\ 
8 & 4 & 0 \\ 
3 & 2 & 7%
\end{array}%
\right] \!,  \notag \\
\mathbf{L}_{3}\left( 4,1,3\right) =\mathbf{\hat{L}}_{3}\mathbf{R}_{3}=\left[ 
\begin{array}{ccc}
5 & 0 & 7 \\ 
6 & 4 & 2 \\ 
1 & 8 & 3%
\end{array}%
\right] \!,\quad \mathbf{L}_{3}\left( 4,1,-3\right) =\mathbf{R}_{3}\mathbf{%
\hat{L}}_{3}^{T}=\left[ 
\begin{array}{ccc}
5 & 6 & 1 \\ 
0 & 4 & 8 \\ 
7 & 2 & 3%
\end{array}%
\right] \!,  \label{L3nat} \\
\mathbf{L}_{3}\left( 4,-1,-3\right) =\mathbf{R}_{3}\mathbf{\hat{L}}_{3}=%
\left[ 
\begin{array}{ccc}
3 & 8 & 1 \\ 
2 & 4 & 6 \\ 
7 & 0 & 5%
\end{array}%
\right] \!,\quad \mathbf{L}_{3}\left( 4,-1,3\right) =\mathbf{\hat{L}}_{3}^{T}%
\mathbf{R}_{3}=\left[ 
\begin{array}{ccc}
3 & 2 & 7 \\ 
8 & 4 & 0 \\ 
1 & 6 & 5%
\end{array}%
\right] \!.  \notag
\end{gather}%
These eight matrices constitute the \emph{phases} \cite{LOLY1} of the magic
square $\mathbf{\hat{L}}_{3}$ which may be regarded as the \emph{fundamental}
natural magic square of\emph{\ }order 3. The classic Lo Shu magic square 
\cite{PICK} is given by $\mathbf{L}_{3}\left( 4,-1,-3\right) +\mathbf{E}%
_{3}. $

On setting%
\begin{equation}
c=v+y+x,  \label{cvy}
\end{equation}%
Lucas' parameterization (\ref{Luc3}) transforms to a general
parameterization given by Frierson \cite{FIER} as%
\begin{align}
\mathbf{G}_{3}\left( x,v,y\right) & =\mathbf{F}_{3}\left( v,y\right) +x%
\mathbf{E}_{3},  \notag \\
\mathbf{F}_{3}\left( v,y\right) & =\left[ 
\begin{array}{ccc}
2v+y & 0 & v+2y \\ 
2y & v+y & 2v \\ 
v & 2v+2y & y%
\end{array}%
\right] =\mathbf{L}_{3}\left( v+y,v,y\right) \!,  \label{F3}
\end{align}%
where $\mathbf{G}_{3}\left( x,v,y\right) $ is essentially equivalent to $%
\mathbf{L}_{3}\left( c,v,y\right) $. However, Frierson takes $x=1$ in
numerical realizations in order to form natural magic squares with elements $%
1,2,\ldots 9$. If $x$ is to be fixed, we prefer $x=0$ for analytical
simplicity. Then, since $v$ and $y$ must be positive, Frierson's
parameterization $\mathbf{F}_{3}\left( v,y\right) $ can generate only two of
the eight phases of the order-3 natural magic square, namely%
\begin{equation}
\mathbf{F}_{3}\left( 3,1\right) =\mathbf{L}_{3}\left( 4,3,1\right) \equiv 
\mathbf{\hat{F}}_{3},\quad \mathbf{F}_{3}\left( 1,3\right) =\mathbf{L}%
_{3}\left( 4,1,3\right) =\mathbf{\hat{F}}_{3}\mathbf{R}_{3},
\end{equation}%
where $\mathbf{\hat{F}}_{3}$ is regarded as fundamental. The other six
phases of $\mathbf{\hat{F}}_{3}$ are not of the form $\mathbf{F}_{3}\left(
v,y\right) $ but can be formed from $\mathbf{G}_{3}\left( x,v,y\right) $.
Furthermore, $\mathbf{L}_{3}$ can be expressed in terms of $\mathbf{F}_{3}$
as%
\begin{equation}
\mathbf{L}_{3}\left( c,v,y\right) =\mathbf{F}_{3}\left( v,y\right) +\left(
c-v-y\right) \mathbf{E}_{3},  \label{L3F3}
\end{equation}%
where $v$ and $y$ need not be positive in this equation. Next, we derive
formulas for the matrices in the JCF and SVD of the matrices $\mathbf{L}%
_{3},\ \mathbf{F}_{3},$ and $\mathbf{E}_{3}$.\smallskip 

\noindent \textbf{JCF.} The Jordan canonical form of a diagonable matrix $%
\mathbf{M}$ is expressed as \cite{MEYE}%
\begin{equation}
\mathbf{M=SDS}^{-1},\quad \mathbf{MS=SD},  \label{Jord}
\end{equation}%
where the eigenvalues are elements of the diagonal matrix $\mathbf{D}$ and
the respective columns of $\mathbf{S}$ are the corresponding eigenvectors.

For $\mathbf{L}_{3}\left( c,v,y\right) ,\ \mathbf{F}_{3}\left( v,y\right) ,$
and $\mathbf{E}_{3},$ we find that%
\begin{align}
\mathbf{S}_{3}\left( v,y\right) & =\left[ 
\begin{array}{ccc}
1 & 1+\Omega \left( v,y\right)  & 1-\Omega \left( v,y\right)  \\ 
1 & -2 & -2 \\ 
1 & 1-\Omega \left( v,y\right)  & 1+\Omega \left( v,y\right) 
\end{array}%
\right] \!,  \notag \\
\mathbf{D}_{L3}\left( c,v,y\right) & =\mathbf{D}_{F3}\left( v,y\right) =%
\limfunc{diag}\left[ \mu _{3},\lambda \left( v,y\right) ,-\lambda \left(
v,y\right) \right] \!,  \label{SJvy} \\
\mathbf{D}_{E3}& =\limfunc{diag}\left[ 3,0,0\right] \!,  \notag
\end{align}%
where%
\begin{equation}
\Omega \left( v,y\right) =3\left( v+y\right) /\lambda \left( v,y\right)
\!,\quad \lambda \left( v,y\right) =\sqrt{3}\sqrt{v^{2}-y^{2}}.  \label{Lam3}
\end{equation}%
Note that $\mathbf{S}_{3}\left( v,y\right) $ and $\lambda \left( v,y\right) $
are independent of $c.$ If $y^{2}>v^{2},$ then there are two imaginary
eigenvalues with imaginary eigenvectors$.$ It is known that regular magic
square matrices have $\pm $ pairs of eigenvalues \cite{MATT,NORD1}%
.\smallskip 

\noindent \textbf{SVD.} The singular value decomposition of a matrix $%
\mathbf{M}$ is expressed as \cite{MEYE}%
\begin{equation}
\mathbf{M=U\Sigma V}^{T},\quad \mathbf{UU}^{T}=\mathbf{I},\quad \mathbf{VV}%
^{T}=\mathbf{I,}  \label{SVD}
\end{equation}%
where $\mathbf{\Sigma }$ is a diagonal matrix with the singular values $%
\sigma _{i}$ as its non-negative elements. The number of non-zero singular
values is equal to the \emph{rank} of the matrix \cite{MEYE}.

For $\mathbf{L}_{3}\left( c,v,y\right) ,\ \mathbf{F}_{3}\left( v,y\right) ,$
and $\mathbf{E}_{3},$ we find that%
\begin{align}
\mathbf{U}_{3}& =\frac{1}{6}\left[ 
\begin{array}{ccc}
2\sqrt{3} & -3\alpha \sqrt{2} & \beta \sqrt{6} \\ 
2\sqrt{3} & 0 & -2\beta \sqrt{6} \\ 
2\sqrt{3} & 3\alpha \sqrt{2} & \beta \sqrt{6}%
\end{array}%
\right] \!,\quad \mathbf{V}_{3}=\frac{1}{6}\left[ 
\begin{array}{ccc}
2\sqrt{3} & -\sqrt{6} & 3\sqrt{2} \\ 
2\sqrt{3} & 2\sqrt{6} & 0 \\ 
2\sqrt{3} & -\sqrt{6} & -3\sqrt{2}%
\end{array}%
\right] \!,  \notag \\
\mathbf{\Sigma }_{L3}\left( c,v,y\right) & =\mathbf{\Sigma }_{F3}\left(
v,y\right) =\limfunc{diag}\left[ \mu _{3},\sigma _{2}\left( v,y\right)
,\sigma _{3}\left( v,y\right) \right] \!,\quad \mathbf{\Sigma }_{E3}=\mathbf{%
D}_{E3},  \label{SVD3}
\end{align}%
where%
\begin{gather}
\sigma _{2}\left( v,y\right) =\alpha \phi \left( v,y\right) =\left\vert \phi
\left( v,y\right) \right\vert \!,\quad \sigma _{3}\left( v,y\right) =\beta
\psi \left( v,y\right) =\left\vert \psi \left( v,y\right) \right\vert \!, 
\notag \\
\phi \left( v,y\right) =\sqrt{3}\left( v+y\right) \!,\quad \psi \left(
v,y\right) =\sqrt{3}\left( v-y\right) \!,  \label{PhiPsi}
\end{gather}%
with $\alpha $ and $\beta $ chosen as $1$ or $-1$ so as to make the singular
values $\sigma _{2}$ and $\sigma _{3}$ non-negative as required when
numerical values are specified for $v$ and $y$. Here, $\mathbf{U}_{3},$ $%
\mathbf{V}_{3},\ \sigma _{2},$ and $\sigma _{3}$ are independent of $c.$%
\smallskip 

\noindent \textbf{Commuting Matrices}. From (\ref{Luc3}), it can be shown
that $\mathbf{L}_{3}\left( c,v,y\right) $ and $\mathbf{L}_{3}\left(
d,s,t\right) $ commute if and only if $vt=ys$ and similarly for $\mathbf{F}%
_{3}\left( v,y\right) $ and $\mathbf{F}_{3}\left( s,t\right) .$ Thus, there
are four pairs of commuting order-3 Lucas matrices among the eight of (\ref%
{L3nat}). The natural Frierson squares $\mathbf{F}_{3}\left( 3,1\right) $
and $\mathbf{F}_{3}\left( 1,3\right) $ do not commute.\smallskip

\noindent \textbf{Matrix Powers.} The matrix powers of $\mathbf{L}_{3}\left(
c,v,y\right) $ can be determined recursively from (\ref{Luc3}) as%
\begin{align}
\mathbf{L}_{3}\left( c,v,y\right) ^{n}& =3^{\frac{n-1}{2}}\left(
v^{2}-y^{2}\right) ^{\frac{n-1}{2}}\left( \mathbf{L}_{3}\left( c,v,y\right)
-c\mathbf{E}_{3}\right) +3^{n-1}c^{n}\mathbf{E}_{3},\text{\quad }n\ \text{%
odd,}  \notag \\
\mathbf{L}_{3}\left( c,v,y\right) ^{n}& =3^{\frac{n}{2}-1}\left(
v^{2}-y^{2}\right) ^{\frac{n}{2}}\left( 3\mathbf{I}_{3}-\mathbf{E}%
_{3}\right) +3^{n-1}c^{n}\mathbf{E}_{3},\text{\quad }n\text{ even.}
\label{L3nth}
\end{align}%
Furthermore, the inverse of $\mathbf{L}_{3}\left( c,v,y\right) $ is given by%
\begin{equation}
\mathbf{L}_{3}\left( c,v,y\right) ^{-1}=\frac{1}{9c}\left(
v^{2}-y^{2}\right) ^{-1}\left( 3c\mathbf{L}_{3}\left( c,v,y\right) +\left(
v^{2}-y^{2}-3c^{2}\right) \mathbf{E}_{3}\right) \!,\text{\quad }v\neq y.
\label{L3inv}
\end{equation}%
These formulas also apply to $\mathbf{F}_{3}\left( v,y\right) ^{n}$ upon
substitution of $c=v+y$ and application of (\ref{F3}).

\section{Frobenius Norm Condition}

The determination of whether a given matrix $\mathbf{M}_{n}$ is natural can
be done by sorting its elements and checking if they are in proper
sequential order $\left( 0,1,\ldots ,n^{2}-1\right) $. Another method is to
check the square of its Frobenius norm $\left\Vert \mathbf{M}_{n}\right\Vert
_{F}^{2}$ (sum of the squares of the matrix elements) to see if it equals
the known sum of elements squared in a natural matrix, i.e. check that%
\begin{equation}
\left\Vert \mathbf{M}_{n}\right\Vert
_{F}^{2}=\tsum\limits_{k=0}^{n^{2}-1}k^{2}=\frac{1}{6}n^{2}\left(
n^{2}-1\right) \left( 2n^{2}-1\right) \!.  \label{Frob1}
\end{equation}%
This is a necessary condition for a matrix to be natural. However, it is not
a sufficient condition as the following counter-example (non-natural) magic
square shows:%
\begin{gather}
\left\Vert \mathbf{M}_{5}\right\Vert _{F}^{2}=\left\Vert \left[ 
\begin{array}{ccccc}
9 & 12 & 15 & 23 & 1 \\ 
15 & 23 & 1 & 9 & 12 \\ 
1 & 9 & 12 & 15 & 23 \\ 
12 & 15 & 23 & 1 & 9 \\ 
23 & 1 & 9 & 12 & 15%
\end{array}%
\right] \right\Vert _{F}^{2}=4900,\quad \limfunc{tr}\left[ \mathbf{M}_{5}%
\right] =60,  \label{M9bad} \\
\frac{1}{6}n^{2}\left( n^{2}-1\right) \left( 2n^{2}-1\right)
|_{n=5}=4900,\quad \mu _{5}=\frac{n}{2}\left( n^{2}-1\right) |_{n=5}=60. 
\notag
\end{gather}%
Thus, the Frobenius norm condition (FNC) is useful mainly as a screening
tool to identify possible natural matrices.

Furthermore, the sum of the squares of the singular values $\sigma _{i}$
equals the Frobenius norm squared \cite{MEYE}, i.e.\footnote{%
This relation also is noted by Cameron, et.\negthinspace\ al \cite{CAM}.}%
\begin{equation}
\left\Vert \mathbf{M}_{n}\right\Vert _{F}^{2}=\tsum\limits_{i=1}^{n}\sigma
_{i}^{2}=\limfunc{tr}\left[ \mathbf{M}_{n}\mathbf{M}_{n}^{T}\right] =\frac{1%
}{6}n^{2}\left( n^{2}-1\right) \left( 2n^{2}-1\right) \!.  \label{Frob2}
\end{equation}%
For a natural magic square matrix, since $\sigma _{1}=\mu _{n},$ (\ref{Frob2}%
) with (\ref{Mun}) becomes%
\begin{equation}
\left\Vert \mathbf{M}_{n}\right\Vert _{F}^{2}=\tsum\limits_{i=2}^{n}\sigma
_{i}^{2}=\frac{1}{12}n^{2}\left( n^{4}-1\right) \!.  \label{Frob3}
\end{equation}%
Thus, the singular values of $\mathbf{M}_{n}$ can be used in the FNC (\ref%
{Frob3}). For example, (\ref{Frob3}) applied to $\mathbf{L}_{3}\left(
c,v,y\right) $ and $\mathbf{F}_{3}\left( v,y\right) $ with their singular
values from (\ref{SVD3}) leads to%
\begin{equation}
v^{2}+y^{2}=10  \label{FrobLF}
\end{equation}%
which is satisfied only for the natural magic squares $\mathbf{L}_{3}\left(
4,\pm 3,\pm 1\right) \!,$ $\mathbf{L}_{3}\left( 4,\pm 1,\pm 3\right) \!,\ 
\mathbf{F}_{3}\left( 3,1\right) \!,$ and $\mathbf{F}_{3}\left( 1,3\right) $
as noted above. For natural $\mathbf{L}_{3}\left( c,v,y\right) $ and $%
\mathbf{F}_{3}\left( v,y\right) ,$ we note that%
\begin{equation}
c=\left\vert v\right\vert +\left\vert y\right\vert =4=\frac{\mu _{3}}{3}%
,\quad v+y=4=\frac{\mu _{3}}{3},  \label{ccvy}
\end{equation}%
respectively.

\section{Compound to Order-9}

\noindent \textbf{Compounding. }As shown in previous work \cite%
{CHAN,EGGE,ROGE,NORD3}, two magic square matrices of order $mn$ can be
formed by compounding two magic square matrices $\mathbf{M}_{m}$ and $%
\mathbf{M}_{n}$ of order $m$ and $n$ as follows:%
\begin{equation}
\mathbf{A}_{mn}=\mathbf{E}_{m}\boldsymbol{\otimes }\mathbf{M}_{n},\quad 
\mathbf{B}_{mn}=\mathbf{M}_{m}\boldsymbol{\otimes }\mathbf{E}_{n},
\label{AB}
\end{equation}%
where the symbol \textquotedblleft $\boldsymbol{\otimes }$%
\textquotedblright\ denotes the Kronecker (tensor) product \cite{MEYE}. It
is shown in \cite{NORD3} that $\mathbf{A}_{mn}$ and $\mathbf{B}_{mn}$
commute and are non-natural magic squares with summation indices $m\mu _{n}$
and $n\mu _{m},$ respectively. If $\mathbf{M}_{m}$ and $\mathbf{M}_{n}$ are
regular, then $\mathbf{A}_{mn}$ and $\mathbf{B}_{mn}$ also are regular.

Application of (\ref{AB}) to the Lucas parameterization $\mathbf{L}_{3}$ of (%
\ref{Luc3}) gives%
\begin{gather}
\mathbf{A}_{9}\left( c,v,y\right) =\mathbf{E}_{3}\boldsymbol{\otimes }%
\mathbf{L}_{3}\left( c,v,y\right)   \notag \\
=\left[ 
\begin{array}{ccccccccc}
v & -v-y & y & v & -v-y & y & v & -v-y & y \\ 
-v+y & 0 & v-y & -v+y & 0 & v-y & -v+y & 0 & v-y \\ 
-y & y+v & -v & -y & y+v & -v & -y & y+v & -v \\ 
v & -v-y & y & v & -v-y & y & v & -v-y & y \\ 
-v+y & 0 & v-y & -v+y & 0 & v-y & -v+y & 0 & v-y \\ 
-y & y+v & -v & -y & y+v & -v & -y & y+v & -v \\ 
v & -v-y & y & v & -v-y & y & v & -v-y & y \\ 
-v+y & 0 & v-y & -v+y & 0 & v-y & -v+y & 0 & v-y \\ 
-y & y+v & -v & -y & y+v & -v & -y & y+v & -v%
\end{array}%
\right] +c\mathbf{E}_{9},  \label{A9}
\end{gather}%
\begin{gather}
\mathbf{B}_{9}\left( d,s,t\right) =\mathbf{L}_{3}\left( d,s,t\right) 
\boldsymbol{\otimes }\mathbf{E}_{3}  \notag \\
=\left[ 
\begin{array}{ccccccccc}
s & s & s & -s-t & -s-t & -s-t & t & t & t \\ 
s & s & s & -s-t & -s-t & -s-t & t & t & t \\ 
s & s & s & -s-t & -s-t & -s-t & t & t & t \\ 
-s+t & -s+t & -s+t & 0 & 0 & 0 & s-t & s-t & s-t \\ 
-s+t & -s+t & -s+t & 0 & 0 & 0 & s-t & s-t & s-t \\ 
-s+t & -s+t & -s+t & 0 & 0 & 0 & s-t & s-t & s-t \\ 
-t & -t & -t & s+t & s+t & s+t & -s & -s & -s \\ 
-t & -t & -t & s+t & s+t & s+t & -s & -s & -s \\ 
-t & -t & -t & s+t & s+t & s+t & -s & -s & -s%
\end{array}%
\right] +d\mathbf{E}_{9}.  \label{B9}
\end{gather}%
It can be verified that $\mathbf{A}\left( c,v,y\right) $ and $\mathbf{B}%
\left( d,s,t\right) $ commute and are non-natural, regular, magic squares
with summation indices $\mu _{A}=9c,\ \mu _{B}=9d.$ Furthermore, they can be
added to produce the order-9 compound Lucas magic square%
\begin{equation}
\mathbf{L}_{9}\left( c,v,y,d,s,t\right) =\mathbf{A}_{9}\left( c,v,y\right) +%
\mathbf{B}_{9}\left( d,s,t\right) \!,\quad \mu _{9}=9\left( c+d\right) \!.
\label{L9}
\end{equation}%
In a similar manner \cite{FIER,LOLY3}, Frierson's order-3 matrix $\mathbf{F}%
_{3}$ of (\ref{F3}) can be compounded to produce the order-9 compound
Frierson magic square%
\begin{equation}
\mathbf{F}_{9}\left( v,y,s,t\right) =\mathbf{E}_{3}\boldsymbol{\otimes }%
\mathbf{F}_{3}\left( v,y\right) +\mathbf{F}_{3}\left( s,t\right) \boldsymbol{%
\otimes }\mathbf{E}_{3}=\mathbf{L}_{9}\left( v+y,v,y,s+t,s,t\right) \!,
\label{FF9}
\end{equation}%
where $v,y,s,t$ must be non-negative. Furthermore, $\mathbf{L}_{9}$ can be
written in terms of $\mathbf{F}_{9}$ as%
\begin{equation}
\mathbf{L}_{9}\left( c,v,y,d,s,t\right) =\mathbf{F}_{9}\left( v,y,s,t\right)
+\left( c+d-v-y-s-t\right) \mathbf{E}_{9},
\end{equation}%
where $v,y,s,t$ may be negative in this equation. In what follows, we
determine numerical values of the parameters for which $\mathbf{L}_{9}\left(
c,v,y,d,s,t\right) $ and $\mathbf{F}_{9}\left( v,y,s,t\right) $ are natural
magic squares.

It is not difficult to show that the phases of $\mathbf{L}_{9}\left(
c,v,y,d,s,t\right) $ are formed in the same manner as those of $\mathbf{L}%
_{3}\left( c,v,y\right) $ in (\ref{L3nat}), e.g.%
\begin{align}
\mathbf{L}_{9}\left( c,v,y,d,s,t\right) \mathbf{R}_{9}& =\mathbf{L}%
_{9}\left( c,y,v,d,t,s\right) \!,  \notag \\
\mathbf{L}_{9}\left( c,v,y,d,s,t\right) ^{T}& =\mathbf{L}_{9}\left(
c,v,-y,d,s,-t\right) \!,  \label{L9ph}
\end{align}%
where the first equation also applies to $\mathbf{F}_{9}\left(
v,y,s,t\right) $ but none of the other phases of $\mathbf{F}_{9}$ are
compound Frierson squares of the form (\ref{FF9}).

Next, we derive the JCF and SVD matrices for $\mathbf{L}_{9}$ and $\mathbf{F}%
_{9}$ based on the work of Nordgren \cite{NORD2,NORD3}. Formulas for JCF and
SVD matrices of compound matrices also are given by Rogers, et.
\negthinspace al \cite{ROGE}. Singular values for compound Frierson squares
are given by Loly and Cameron \cite{LOLY3}.\smallskip

\noindent \textbf{JCF.} We recall from \cite{NORD3} that the JCF matrices of 
$\mathbf{A}_{mn},$ $\mathbf{B}_{mn},$ and $\mathbf{E}_{mn}$ are given by%
\begin{equation}
\mathbf{A}_{mn}=\mathbf{S}_{mn}\mathbf{D}_{Amn}\mathbf{S}_{mn}^{-1},\quad 
\mathbf{B}_{mn}=\mathbf{S}_{mn}\mathbf{D}_{Bmn}\mathbf{S}_{mn}^{-1},\quad 
\mathbf{E}_{mn}=\mathbf{S}_{mn}\mathbf{D}_{Emn}\mathbf{S}_{mn}^{-1},
\label{JordAB}
\end{equation}%
where%
\begin{gather}
\mathbf{S}_{mn}=\mathbf{S}_{m}\boldsymbol{\otimes }\mathbf{S}_{n},\quad 
\mathbf{D}_{Amn}=\mathbf{D}_{Em}\boldsymbol{\otimes }\mathbf{D}_{Mn},  \notag
\\
\mathbf{D}_{Bmn}=\mathbf{D}_{Mm}\boldsymbol{\otimes }\mathbf{D}_{En},\quad 
\mathbf{D}_{Emn}=\limfunc{diag}\left[ mn,0,0,\ldots ,0\right] \!.
\label{SJAB}
\end{gather}%
Application of (\ref{SJAB}) to $\mathbf{A}_{9}\left( c,v,y\right) ,\ \mathbf{%
B}_{9}\left( d,s,t\right) ,$ $\mathbf{L}_{9}\left( c,v,y,d,s,t\right) ,$ and 
$\mathbf{E}_{9}$ gives%
\begin{align}
\mathbf{S}_{9}& =\mathbf{S}_{3}\left( v,y\right) \boldsymbol{\otimes }%
\mathbf{S}_{3}\left( s,t\right) \!,\quad \mathbf{D}_{A9}=\mathbf{D}_{E3}%
\boldsymbol{\otimes }\mathbf{D}_{L3}\left( c,v,y\right) \!,  \notag \\
\mathbf{D}_{B9}& =\mathbf{D}_{L3}\left( d,s,t\right) \boldsymbol{\otimes }%
\mathbf{D}_{E3},\quad \mathbf{D}_{L9}=\mathbf{D}_{A9}+\mathbf{D}_{B9},\quad 
\mathbf{D}_{E9}=\limfunc{diag}\left[ 9,0,0,\ldots ,0\right] \!,  \label{S9J9}
\end{align}%
and similarly for $\mathbf{F}_{9}\left( v,y,s,t\right) .$ From (\ref{S9J9}),
after some manipulation using (\ref{SJvy}), we find the following JCF
matrices for $\mathbf{L}_{9},\ \mathbf{F}_{9},$ and $\mathbf{E}_{9}$:%
\begin{equation}
\mathbf{S}_{9}\left( v,y,s,t\right) =\left[ 
\begin{array}{ccccccccc}
1 & 1+\Omega \left( v,y\right)  & 1-\Omega \left( v,y\right)  & 1+\Omega
\left( s,t\right)  & 1-\Omega \left( s,t\right)  & 1 & 0 & -1 & 0 \\ 
1 & -2 & -2 & 1+\Omega \left( s,t\right)  & 1-\Omega \left( s,t\right)  & -1
& -1 & 0 & 0 \\ 
1 & 1-\Omega \left( v,y\right)  & 1+\Omega \left( v,y\right)  & 1+\Omega
\left( s,t\right)  & 1-\Omega \left( s,t\right)  & 0 & 1 & 1 & 0 \\ 
1 & 1+\Omega \left( v,y\right)  & 1-\Omega \left( v,y\right)  & -2 & -2 & -1
& 1 & 0 & -1 \\ 
1 & -2 & -2 & -2 & -2 & 0 & 0 & 0 & 1 \\ 
1 & 1-\Omega \left( v,y\right)  & 1+\Omega \left( v,y\right)  & -2 & -2 & 1
& -1 & 0 & 0 \\ 
1 & 1+\Omega \left( v,y\right)  & 1-\Omega \left( v,y\right)  & 1-\Omega
\left( s,t\right)  & 1+\Omega \left( s,t\right)  & 0 & -1 & 1 & 1 \\ 
1 & -2 & -2 & 1-\Omega \left( s,t\right)  & 1+\Omega \left( s,t\right)  & 1
& 1 & 0 & -1 \\ 
1 & 1-\Omega \left( v,y\right)  & 1+\Omega \left( v,y\right)  & 1-\Omega
\left( s,t\right)  & 1+\Omega \left( s,t\right)  & -1 & 0 & -1 & 0%
\end{array}%
\right] \!,  \label{S9}
\end{equation}%
\begin{align}
\mathbf{D}_{L9}\left( c,v,y,d,s,t\right) & =\mathbf{D}_{F9}\left(
v,y,s,t\right) =\limfunc{diag}\left[ \mu _{9},3\lambda \left( v,y\right)
\!,-3\lambda \left( v,y\right) \!,3\lambda \left( s,t\right) \!,-3\lambda
\left( s,t\right) \!,0,0,0,0\right] \!,  \notag \\
\mathbf{D}_{E9}& =\limfunc{diag}\left[ 9,0,0,0,0,0,0,0,0\right] \!,
\label{D9}
\end{align}%
where $\Omega $ and $\lambda $ are given by (\ref{Lam3}).\smallskip 

\noindent \textbf{SVD.} We recall from \cite{NORD3} that the SVD of $\mathbf{%
A}_{mn},$ $\mathbf{B}_{mn},$ and $\mathbf{E}_{mn}$ are given by 
\begin{equation}
\mathbf{A}_{mn}=\mathbf{U}_{mn}\mathbf{\Sigma }_{Amn}\mathbf{V}%
_{mn}^{T},\quad \mathbf{B}_{mn}=\mathbf{U}_{mn}\mathbf{\Sigma }_{Bmn}\mathbf{%
V}_{mn}^{T},\quad \mathbf{E}_{mn}=\mathbf{U}_{mn}\mathbf{\Sigma }_{Emn}%
\mathbf{V}_{mn}^{T},  \label{ABsvd}
\end{equation}%
where%
\begin{gather*}
\mathbf{U}_{mn}=\mathbf{U}_{m}\boldsymbol{\otimes }\mathbf{U}_{n},\quad 
\mathbf{V}_{mn}=\mathbf{V}_{m}\boldsymbol{\otimes }\mathbf{V}_{n},\quad 
\mathbf{\Sigma }_{Emn}=\mathbf{D}_{Emn}, \\
\mathbf{\Sigma }_{Amn}=\mathbf{\Sigma }_{Em}\boldsymbol{\otimes }\mathbf{%
\Sigma }_{Mn},\quad \mathbf{\Sigma }_{Bmn}=\mathbf{\Sigma }_{Mm}\boldsymbol{%
\otimes }\mathbf{\Sigma }_{En}
\end{gather*}%
from which we find the following SVD matrices\footnote{%
The columns of these SVD matrices are rearranged for clarity.} for $\mathbf{L%
}_{9},\ \mathbf{F}_{9}$ and $\mathbf{E}_{9}$:%
\begin{equation}
\mathbf{U}_{9}=\mathbf{U}_{3}\boldsymbol{\otimes }\mathbf{U}_{3}=\frac{1}{6}%
\left[ 
\begin{array}{ccccccccc}
2 & -\alpha \sqrt{6} & \beta \sqrt{2} & -\gamma \sqrt{6} & \eta \sqrt{2} & -%
\sqrt{3} & 1 & -\sqrt{3} & 3 \\ 
2 & 0 & -2\beta \sqrt{2} & -\gamma \sqrt{6} & \eta \sqrt{2} & 0 & -2 & 2%
\sqrt{3} & 0 \\ 
2 & \alpha \sqrt{6} & \beta \sqrt{2} & -\gamma \sqrt{6} & \eta \sqrt{2} & 
\sqrt{3} & 1 & -\sqrt{3} & -3 \\ 
2 & -\alpha \sqrt{6} & \beta \sqrt{2} & 0 & -2\eta \sqrt{2} & 2\sqrt{3} & -2
& 0 & 0 \\ 
2 & 0 & -2\beta \sqrt{2} & 0 & -2\eta \sqrt{2} & 0 & 4 & 0 & 0 \\ 
2 & \alpha \sqrt{6} & \beta \sqrt{2} & 0 & -2\eta \sqrt{2} & -2\sqrt{3} & -2
& 0 & 0 \\ 
2 & -\alpha \sqrt{6} & \beta \sqrt{2} & \gamma \sqrt{6} & \eta \sqrt{2} & -%
\sqrt{3} & 1 & \sqrt{3} & -3 \\ 
2 & 0 & -2\beta \sqrt{2} & \gamma \sqrt{6} & \eta \sqrt{2} & 0 & -2 & -2%
\sqrt{3} & 0 \\ 
2 & \alpha \sqrt{6} & \beta \sqrt{2} & \gamma \sqrt{6} & \eta \sqrt{2} & 
\sqrt{3} & 1 & \sqrt{3} & 3%
\end{array}%
\right] \!,  \label{U9}
\end{equation}%
\begin{equation}
\mathbf{V}_{9}=\mathbf{V}_{3}\boldsymbol{\otimes }\mathbf{V}_{3}=\frac{1}{6}%
\left[ 
\begin{array}{ccccccccc}
2 & -\sqrt{2} & \sqrt{6} & -\sqrt{2} & \sqrt{6} & -\sqrt{3} & 3 & -\sqrt{3}
& 1 \\ 
2 & 2\sqrt{2} & 0 & -\sqrt{2} & \sqrt{6} & 2\sqrt{3} & 0 & 0 & -2 \\ 
2 & -\sqrt{2} & -\sqrt{6} & -\sqrt{2} & \sqrt{6} & -\sqrt{3} & -3 & \sqrt{3}
& 1 \\ 
2 & -\sqrt{2} & \sqrt{6} & 2\sqrt{2} & 0 & 0 & 0 & 2\sqrt{3} & -2 \\ 
2 & 2\sqrt{2} & 0 & 2\sqrt{2} & 0 & 0 & 0 & 0 & 4 \\ 
2 & -\sqrt{2} & -\sqrt{6} & 2\sqrt{2} & 0 & 0 & 0 & -2\sqrt{3} & -2 \\ 
2 & -\sqrt{2} & \sqrt{6} & -\sqrt{2} & -\sqrt{6} & \sqrt{3} & -3 & -\sqrt{3}
& 1 \\ 
2 & 2\sqrt{2} & 0 & -\sqrt{2} & -\sqrt{6} & -2\sqrt{3} & 0 & 0 & -2 \\ 
2 & -\sqrt{2} & -\sqrt{6} & -\sqrt{2} & -\sqrt{6} & \sqrt{3} & 3 & \sqrt{3}
& 1%
\end{array}%
\right] \!,  \label{V9}
\end{equation}

$\vspace{-0.1in}$%
\begin{align}
\mathbf{\Sigma }_{L9}\left( c,v,y,d,s,t\right) & =\mathbf{\Sigma }%
_{A9}\left( c,v,y\right) +\mathbf{\Sigma }_{B9}\left( d,s,t\right) =\mathbf{%
\Sigma }_{E3}\mathbf{\boldsymbol{\otimes }\Sigma }_{L3}\left( c,v,y\right) +%
\mathbf{\Sigma }_{L3}\left( d,s,t\right) \mathbf{\boldsymbol{\otimes }\Sigma 
}_{E3}  \notag \\
& =\limfunc{diag}\left[ \mu _{9},3\alpha \phi \left( v,y\right) \!,3\beta
\psi \left( v,y\right) \!,3\gamma \phi \left( s,t\right) \!,3\eta \psi
\left( s,t\right) \!,0,0,0,0\right]   \notag \\
& =\limfunc{diag}\left[ \mu _{9},3\left\vert \phi \left( v,y\right)
\right\vert \!,3\left\vert \psi \left( v,y\right) \right\vert \!,3\left\vert
\phi \left( s,t\right) \right\vert \!,3\left\vert \psi \left( s,t\right)
\right\vert \!,0,0,0,0\right] \!,  \notag \\
\mathbf{\Sigma }_{F9}\left( v,y,s,t\right) & =\mathbf{\Sigma }_{L9}\left(
c,v,y,d,s,t\right) ,\quad \mathbf{\Sigma }_{E9}=\mathbf{D}_{E9},
\label{Sig9}
\end{align}%
where and $\phi $ and $\psi $ are given by (\ref{PhiPsi}) and $\alpha ,\beta
,\gamma ,\eta $ are chosen as $1$ or $-1$ so as to make the singular values
positive as required when numerical values are specified for $v,y,s,t$. The
singular values in $\mathbf{\Sigma }_{F9}$ agree with those given by Loly
and Cameron \cite{LOLY3} and $\mathbf{\Sigma }_{F9}$ is the same as $\mathbf{%
\Sigma }_{L9}$ since the $\pm $ signs on $v,y,s,t$ have no effect. We note
that $\mathbf{L}_{9}$ and $\mathbf{F}_{9}$ are rank 5 matrices.\smallskip 

\noindent \textbf{Natural Squares}. The Frobenius norm condition (FNC) of (%
\ref{Frob3}) applied to the singular values $\mathbf{\Sigma }_{L9}$ and $%
\mathbf{\Sigma }_{F9}$ of (\ref{Sig9}) leads to%
\begin{equation}
v^{2}+y^{2}+s^{2}+t^{2}=820  \label{vyst}
\end{equation}%
which is satisfied by $v,y,s,t$ taking the values $\pm 1,\pm 3,\pm 9,\pm 27$
in any order. It can be verified that these values produce natural magic
squares $\mathbf{L}_{9}\left( c,v,y,d,s,t\right) $ with 
\begin{equation}
c=\left\vert v\right\vert +\left\vert y\right\vert \!,\quad d=\left\vert
s\right\vert +\left\vert t\right\vert \!,\quad c+d=\frac{\mu _{9}}{9}=40.
\label{cdmu}
\end{equation}%
A direct search shows that no other distinct integer solutions of (\ref{vyst}%
) exist. Therefore, no other order-9 natural compound Lucas squares can be
formed. Similarly, when $v,y,s,t$ take the values $1,3,9,27$ in any order,
the Frierson square $\mathbf{F}_{9}\left( v,y,s,t\right) $ also is natural.

The naturalness of $\mathbf{L}_{9}$ also can be established by noting that $%
\mathbf{A}_{9}\left( c,v,y\right) $ and $\mathbf{B}_{9}\left( d,s,t\right) $
of (\ref{A9}) and (\ref{B9}) are mutually orthogonal in form and their sum
may produce a $\mathbf{L}_{9}\left( c,v,y,d,s,t\right) $ without duplicate
elements for certain values of $v,y,s,t.$ In order to determine if there are
any duplicate elements, we observe that all elements of a $\mathbf{L}_{9}$
that satisfies the FNC are of the form%
\begin{equation}
\alpha _{0}+3\alpha _{1}+9\alpha _{2}+27\alpha _{3}+40,\quad \alpha
_{i}=-1,0,\text{ or }1.
\end{equation}%
Thus, two elements are identical if $\alpha _{i}$ and $\beta _{i}$ can be
found such that 
\begin{equation}
\alpha _{0}+3\alpha _{1}+9\alpha _{2}+27\alpha _{3}=\beta _{0}+3\beta
_{1}+9\beta _{2}+27\beta _{3},\quad \alpha _{i},\beta _{i}=-1,0,\text{ or }1,
\end{equation}%
i.e.%
\begin{equation}
\Lambda \equiv \eta _{0}+3\eta _{1}+9\eta _{2}+27\eta _{3}=0,\quad \eta
_{i}=\alpha _{i}-\beta _{i}=-2,-1,0,1,\text{ or }2.
\end{equation}%
The critical cases are 
\begin{align}
\eta _{0}& =\eta _{1}=\eta _{2}=-2,\ \eta _{3}=1,\;\Rightarrow \;\Lambda =1,
\notag \\
\eta _{0}& =\eta _{1}=-2,\ \eta _{2}=1,\ \eta _{3}=0,\;\Rightarrow \;\Lambda
=1, \\
\eta _{0}& =-2,\text{\ }\eta _{1}=1,\ \eta _{2}=\eta _{3}=0,\;\Rightarrow
\;\Lambda =1,  \notag
\end{align}%
hence, there are no duplicate elements in a $\mathbf{L}_{9}$ that satisfies
the FNC. Therefore, $\mathbf{L}_{9}\left( c,v,y,d,s,t\right) $ is natural
for $v,y,s,t$ taking the values $\pm 1,\pm 3,\pm 9,\pm 27$ in any order with 
$c$ and $d$ from (\ref{cdmu}). Similarly, $\mathbf{F}_{9}\left(
v,y,s,t\right) $ is natural for $v,y,s,t$ taking the values $1,3,9,27$ in
any order.

There are $4!$ orderings of $1,3,9,27$. For Frierson squares, a phase square
results only from the first of (\ref{L9ph}), namely%
\begin{equation}
\mathbf{F}_{9}\left( y,v,t,s\right) =\mathbf{F}_{9}\left( v,y,s,t\right) 
\mathbf{R}_{9},  \label{F9phase}
\end{equation}%
giving $4!/2=12$ as the number of fundamental (phases excluded) natural
compound Frierson magic squares of order 9. However, Frierson \cite{FIER}
gives only six such squares, namely those equivalent to%
\begin{gather}
\mathbf{F}_{9A}=\mathbf{F}_{9}\left( 3,1,27,9\right) \!,\quad \mathbf{F}%
_{9D}=\mathbf{F}_{9}\left( 27,9,3,1\right) \!,  \notag \\
\mathbf{F}_{9B}=\mathbf{F}_{9}\left( 27,1,9,3\right) \!,\quad \mathbf{F}%
_{9E}=\mathbf{F}_{9}\left( 9,3,27,1\right) \!,  \label{Fier9} \\
\mathbf{F}_{9C}=\mathbf{F}_{9}\left( 9,1,27,3\right) \!,\quad \mathbf{F}%
_{9F}=\mathbf{F}_{9}\left( 27,3,9,1\right) \!,  \notag
\end{gather}%
the other six being%
\begin{gather}
\mathbf{F}_{9G}=\mathbf{F}_{9}\left( 3,1,9,27\right) \!,\quad \mathbf{F}%
_{9J}=\mathbf{F}_{9}\left( 9,27,3,1\right) \!,  \notag \\
\mathbf{F}_{9H}=\mathbf{F}_{9}\left( 27,1,3,9\right) \!,\quad \mathbf{F}%
_{9K}=\mathbf{F}_{9}\left( 3,9,27,1\right) \!,  \label{Fier9+} \\
\mathbf{F}_{9I}=\mathbf{F}_{9}\left( 9,1,3,27\right) \!,\quad \mathbf{F}%
_{9L}=\mathbf{F}_{9}\left( 3,27,9,1\right)  \notag
\end{gather}%
which are not mere phases of Frierson's six. However, Frierson and others 
\cite{BELL,LOLY3} do not consider those of (\ref{Fier9+}) to be fundamental.

For compound Lucas squares there are $2^{4}$ possible $\pm $ sign
arrangements for each of the $4!$ orderings of $1,3,9,27$ for a total of $%
384 $ natural Lucas squares. Since for each fundamental Lucas square there
are $7 $ phase squares included in the $384$, the number of fundamental
order-9, natural Lucas squares is $384/8=48.$ A detailed examination of the $%
384$ Lucas squares confirms that $48$ are fundamental and the remainder are
phases of them as expected. We find that $48$ fundamental Lucas squares can
be formed from $\mathbf{L}_{9}\left( c,v,y,d,s,t\right) $, $\mathbf{L}%
_{9}\left( c,v,y,d,s,-t\right) $, $\mathbf{L}_{9}\left( c,v,y,d,-s,-t\right) 
$, and $\mathbf{L}_{9}\left( c,v,-y,d,-s,-t\right) $ using the $v,y,s,t$
values of the $12$ fundamental Frierson squares (\ref{Fier9}) and (\ref%
{Fier9+}) together with (\ref{cdmu}).

Nonzero eigenvalues and singular values of the fundamental $\mathbf{L}_{9}$%
's and $\mathbf{F}_{9}$'s are given in Table 1.

\begin{gather*}
\begin{tabular}{||c||c|c||c|c|c|c||}
\hline\hline
$%
\begin{array}{c}
v,y,s,t \\ 
\text{from }\mathbf{F}_{9}%
\end{array}%
$ & $\left\vert \mathbf{\lambda }_{1}\right\vert $ & $\left\vert \mathbf{%
\lambda }_{2}\right\vert $ & $\sigma _{2}/\sqrt{3}$ & $\sigma _{3}/\sqrt{3}$
& $\sigma _{4}/\sqrt{3}$ & $\sigma _{5}/\sqrt{3}$ \\ \hline
$A,\ G$ & $6\sqrt{6}$ & $54\sqrt{6}$ & $12$ & $6$ & $108$ & $54$ \\ 
$D,\ J$ & $54\sqrt{6}$ & $6\sqrt{6}$ & $108$ & $54$ & $12$ & $6$ \\ 
$B,\ H$ & $\ 6\sqrt{546}\ $ & $18\sqrt{6}$ & $84$ & $78$ & $36$ & $18$ \\ 
$E,\ K$ & $18\sqrt{6}$ & $\ 6\sqrt{546}\ $ & $36$ & $18$ & $84$ & $78$ \\ 
$C,\ I$ & $12\sqrt{15}$ & $36\sqrt{15}$ & $30$ & $24$ & $90$ & $72$ \\ 
$F,\ L$ & $36\sqrt{15}$ & $12\sqrt{15}$ & $90$ & $72$ & $30$ & $24$ \\ 
\hline\hline
\end{tabular}
\\
\text{Table 1 - EV's and SV's of order-9 Lucas and Frierson matrices.}
\end{gather*}

\noindent \textbf{Commuting Matrices}. From (\ref{L9}), (\ref{A9}), and (\ref%
{B9}), it can be shown that 
\begin{align}
\left[ \mathbf{L}_{9}\left( c,v,y,d,s,t\right) ,\mathbf{L}_{9}\left(
d,s,t,c,v,y\right) \right] & =\mathbf{0}\;\text{iff}\;vt=ys,  \notag \\
\left[ \mathbf{L}_{9}\left( c,v,y,d,s,t\right) ,\mathbf{L}_{9}\left(
d,t,s,c,y,v\right) \right] & =\mathbf{0}\;\text{iff}\;vs=yt,  \label{L9L9}
\end{align}%
where $\left[ \mathbf{A,B}\right] =\mathbf{AB-BA}$ is the \emph{commutator}
of $\mathbf{A}$ and $\mathbf{B}$. Similar results hold for $\mathbf{F}%
_{9}\left( v,y,s,t\right) $ and thus, only the following eight pairs of $%
\mathbf{F}_{9}$'s from (\ref{Fier9}) and (\ref{Fier9+}) commute:%
\begin{gather}
\left[ \mathbf{F}_{9A},\mathbf{F}_{9D}\right] =\left[ \mathbf{F}_{9C},%
\mathbf{F}_{9F}\right] =\left[ \mathbf{F}_{9A}\mathbf{R}_{9},\mathbf{F}_{9D}%
\mathbf{R}_{9}\right] =\left[ \mathbf{F}_{9C}\mathbf{R}_{9},\mathbf{F}_{9F}%
\mathbf{R}_{9}\right] =\mathbf{0},  \notag \\
\left[ \mathbf{F}_{9G},\mathbf{F}_{9J}\mathbf{R}_{9}\right] =\left[ \mathbf{F%
}_{9G}\mathbf{R}_{9},\mathbf{F}_{9J}\right] =\left[ \mathbf{F}_{9I},\mathbf{F%
}_{9L}\mathbf{R}_{9}\right] =\left[ \mathbf{F}_{9I}\mathbf{R}_{9},\mathbf{F}%
_{9L}\right] =\mathbf{0},  \label{F9comm}
\end{gather}%
where only the first two commutators involve fundamental $\mathbf{F}_{9}$'s.
The conditions in (\ref{L9L9}) can be verified by observing that the
eigenvector matrix (\ref{S9}) of the commuting $\mathbf{L}_{9}$'s are
identical under them and the $\mathbf{L}_{9}$'s are diagonable. Furthermore,
if any regular (associative) magic squares $\mathbf{A}$ and $\mathbf{B}$
commute, then all like phases of $\mathbf{A}$ and $\mathbf{B}$ commute,
including $\mathbf{AR}$ and $\mathbf{BR}$ as in (\ref{F9comm}). For the four
forms of fundamental Lucas squares noted above, only $\mathbf{L}_{9}\left(
c,v,y,d,s,t\right) $ and $\mathbf{L}_{9}\left( c,v,y,d,-s,-t\right) $ (with $%
v,y,s,t$ $\geq 0$) and their eight phases give rise to 64 pairs of commuting
matrices according to (\ref{L9L9}) and (\ref{F9comm}).\smallskip 

\noindent \textbf{Matrix Powers.} The matrix powers of $\mathbf{L}_{9}\left(
c,v,y,d,s,t\right) $ can be determined recursively from (\ref{L9}) and (\ref%
{L3nth}) as 
\begin{gather}
\mathbf{L}_{9}\left( c,v,y,d,s,t\right) ^{n}=\left( 3\sqrt{3}\right)
^{n-1}\left( v^{2}-y^{2}\right) ^{\frac{n-1}{2}}\left( \mathbf{A}_{9}\left(
c,v,y\right) -c\mathbf{E}_{9}\right)   \notag \\
+\left( 3\sqrt{3}\right) ^{n-1}\left( s^{2}-t^{2}\right) ^{\frac{n-1}{2}%
}\left( \mathbf{B}_{9}\left( d,s,t\right) -d\mathbf{E}_{9}\right)
+9^{n-1}\left( c+d\right) ^{n}\mathbf{E}_{9},\text{\quad }n\text{ odd,} 
\notag \\
\vspace{-0.2in} \\
\mathbf{L}_{9}\left( c,v,y,d,s,t\right) ^{n}=3^{\frac{3n}{2}-2}\left(
v^{2}-y^{2}\right) ^{\frac{n}{2}}\left( 3\mathbf{E}_{3}\mathbf{\boldsymbol{%
\otimes }I}_{3}-\mathbf{E}_{9}\right)   \notag \\
+\ 3^{\frac{3n}{2}-2}\left( s^{2}-t^{2}\right) ^{\frac{n}{2}}\left( 3\mathbf{%
I}_{3}\mathbf{\boldsymbol{\otimes }E}_{3}-\mathbf{E}_{9}\right)
+9^{n-1}\left( c+d\right) ^{n}\mathbf{E}_{9},\text{\quad }n\text{ even.} 
\notag
\end{gather}%
These formulas also apply to $\mathbf{F}_{9}\left( v,y,s,t\right) ^{n}$ upon
substitution of $c=v+y$ and $d=s+t.$ Since $\mathbf{L}_{9}$ and $\mathbf{F}%
_{9}$ are singular, they have no inverse.

\section{Sequential Compounding}

\noindent \textbf{Compounding. }In this section only, $\mathbf{M}_{\ell }$
denotes a matrix of order $n=3^{\ell },$ where $\ell $ is called the \emph{%
level} of $\mathbf{M}_{\ell }$ as in \cite{LOLY3}. On applying (\ref{AB}) in
a sequential manner, $\mathbf{A}_{\ell }$ and $\mathbf{B}_{\ell }$ can be
formed as%
\begin{equation}
\mathbf{A}_{\ell }=\mathbf{E}_{1}\boldsymbol{\otimes }\mathbf{L}_{\ell
-1},\quad \mathbf{B}_{\ell }=\mathbf{L}_{1}\boldsymbol{\otimes }\mathbf{E}%
_{\ell -1},\quad \ell =2,3,\ldots \ ,  \label{AB3n}
\end{equation}%
where $\mathbf{L}_{1}$ is given by (\ref{Luc3}) and a compound Lucas square
is formed as\footnote{%
This compounding construction appears to be equivalent to that of Loly and
Cameron \cite{LOLY3} for Frierson squares, although their symbol
\textquotedblleft $\boldsymbol{\otimes }$\textquotedblright\ does not denote
the standard Kronecker product.}%
\begin{equation}
\mathbf{L}_{\ell }=\mathbf{A}_{\ell }+\mathbf{B}_{\ell }.  \label{L3n}
\end{equation}%
New parameter names for the parameters in $\mathbf{L}_{1}\left( c,v,y\right) 
$ must be used in (\ref{AB3n}) at each step, i.e.%
\begin{equation}
\mathbf{L}_{\ell +1}=\mathbf{E}_{1}\boldsymbol{\otimes }\mathbf{L}_{\ell
}\left( c_{1},v_{1},y_{1},c_{2},v_{2},y_{2},\ldots ,c_{\ell },v_{\ell
},y_{\ell }\right) +\mathbf{L}_{1}\left( c_{\ell +1},v_{\ell +1},y_{\ell
+1}\right) \boldsymbol{\otimes }\mathbf{E}_{\ell }.  \label{Ll+1}
\end{equation}%
It is not difficult to show by induction that $\mathbf{A}_{\ell }$ and $%
\mathbf{B}_{\ell }$ commute, are magic and regular, and that the phases of $%
\mathbf{L}_{\ell }\left( c_{i},v_{i},y_{i}\right) ,$ $\left( i=1,2,\ldots
,\ell \right) ,$ are formed in the same manner as those of $\mathbf{L}%
_{3}\left( c,v,y\right) $ in (\ref{L3nat}), e.g.%
\begin{equation}
\mathbf{L}_{\ell }\left( c_{i},v_{i},y_{i}\right) \mathbf{R}_{\ell }=\mathbf{%
L}_{\ell }\left( c_{i},y_{i},v_{i}\right) \!,\quad \mathbf{L}_{\ell
}^{T}\left( c_{i},v_{i},y_{i}\right) =\mathbf{L}_{\ell }\left(
c_{i},-v_{i},-y_{i}\right) \!.  \label{L2ph}
\end{equation}%
A Frierson square $\mathbf{F}_{\ell }\left( v_{1},y_{1},v_{2},y_{2},\ldots
,v_{\ell },y_{\ell }\right) $ can be constructed by compounding in the same
manner as $\mathbf{L}_{\ell }$ above, but only the first form of (\ref{L2ph}%
) applies to it. Furthermore, $\mathbf{F}_{\ell }\left( v_{i},y_{i}\right) $
and $\mathbf{L}_{\ell }\left( c_{i},v_{i},y_{i}\right) $ are related by 
\begin{align}
\mathbf{F}_{\ell }\left( v_{i},y_{i}\right) & =\mathbf{L}_{\ell }\left(
v_{i}+y_{i},v_{i},y_{i}\right) \!,\quad v_{i}\geq 0,\ y_{i}\geq 0,\quad
i=1,2,\ldots ,\ell ,  \notag \\
\mathbf{L}_{\ell }\left( c_{i},v_{i},y_{i}\right) & =\mathbf{F}_{\ell
}\left( v_{i},y_{i}\right) +\tsum\limits_{i=1}^{\ell }\left(
c_{i}-v_{i}-y_{i}\right) \mathbf{E}_{\ell },\quad v_{i},\ y_{i}\text{
unrestricted.}  \label{FlLl}
\end{align}%
In what follows we determine values for $c_{i},v_{i},y_{i}$ such that $%
\mathbf{L}_{\ell }\left( c_{i},v_{i},y_{i}\right) $ and $\mathbf{F}_{\ell
}\left( v_{i},y_{i}\right) $ are natural magic squares. First, we derive
their JCF and SVD.\smallskip 

\noindent \textbf{JCF.} According to (\ref{SJAB}) with (\ref{AB3n}), the
eigenvector and eigenvalue matrices of $\mathbf{A}_{\ell },$ $\mathbf{B}%
_{\ell },\ \mathbf{L}_{\ell },$ and $\mathbf{E}_{\ell }$ are given by 
\begin{equation}
\left. 
\begin{array}{c}
\mathbf{S}_{\ell }=\mathbf{S}_{1}\boldsymbol{\otimes }\mathbf{S}_{\ell
-1},\quad \mathbf{D}_{A\ell }=\mathbf{D}_{E1}\boldsymbol{\otimes }\mathbf{D}%
_{L,\ell -1},\quad \mathbf{D}_{B\ell }=\mathbf{D}_{L1}\boldsymbol{\otimes }%
\mathbf{D}_{E,\ell -1},\vspace{0.05in} \\ 
\mathbf{D}_{L\ell }=\mathbf{D}_{A\ell }+\mathbf{D}_{B\ell },\quad \mathbf{D}%
_{E\ell }=\limfunc{diag}\left[ 3^{\ell },0,0,\ldots ,0\right]%
\end{array}%
\right\} \quad \ell =2,3,\ldots \ .  \label{Jordn1}
\end{equation}%
This can be verified sequentially as follows:%
\begin{align}
\mathbf{L}_{\ell }& =\mathbf{S}_{\ell }\mathbf{D}_{L\ell }\mathbf{S}_{\ell
}^{-1}=\left( \mathbf{S}_{1}\boldsymbol{\otimes }\mathbf{S}_{\ell -1}\right)
\left( \mathbf{D}_{E1}\boldsymbol{\otimes }\mathbf{D}_{L,\ell -1}+\mathbf{D}%
_{L1}\boldsymbol{\otimes }\mathbf{D}_{E,\ell -1}\right) \left( \mathbf{S}%
_{1}^{-1}\boldsymbol{\otimes }\mathbf{S}_{\ell -1}^{-1}\right)  \notag \\
& =\left( \mathbf{S}_{1}\mathbf{D}_{L1}\mathbf{S}_{1}^{-1}\right) 
\boldsymbol{\otimes }\left( \mathbf{S}_{\ell -1}\mathbf{D}_{L,\ell -1}%
\mathbf{S}_{\ell -1}^{-1}\right) +\left( \mathbf{S}_{1}\mathbf{D}_{L1}%
\mathbf{S}_{1}^{-1}\right) \boldsymbol{\otimes }\left( \mathbf{S}_{\ell -1}%
\mathbf{D}_{E,\ell -1}\mathbf{S}_{\ell -1}^{-1}\right)  \label{Jordl1} \\
& =\mathbf{E}_{1}\boldsymbol{\otimes }\mathbf{L}_{\ell -1}+\mathbf{L}_{1}%
\boldsymbol{\otimes }\mathbf{E}_{\ell -1}=\mathbf{A}_{\ell }+\mathbf{B}%
_{\ell }=\mathbf{L}_{\ell },\quad \ell =2,3,\ldots \ .  \notag
\end{align}%
From (\ref{Jordn1}), starting with (\ref{SJvy}) and (\ref{Sig9}), the
nonzero eigenvalues of $\mathbf{L}_{\ell }\left( c_{i},v_{i},y_{i}\right) $
are found to be%
\begin{equation}
\mu _{n},\ \pm 3^{\ell -1}\lambda \left( v_{1},y_{1}\right) \!,\ \pm 3^{\ell
-1}\lambda \left( v_{2},y_{2}\right) \!,\ldots ,\ \pm 3^{\ell -1}\lambda
\left( v_{\ell },y_{\ell }\right) \!,\quad \ell =1,2,\ldots \ ,
\end{equation}%
with $n=3\ell $ and $\lambda \left( v_{i},y_{i}\right) $ from (\ref{Lam3}).
These eigenvalues also apply to compound Frierson matrices.\smallskip

\noindent \textbf{SVD.} According to (\ref{ABsvd}) with (\ref{AB3n}), the
SVD matrices of $\mathbf{A}_{\ell },$ $\mathbf{B}_{\ell },\ \mathbf{L}_{\ell
},$ and $\mathbf{E}_{\ell }$ are given by%
\begin{equation}
\left. 
\begin{array}{c}
\mathbf{U}_{\ell }=\mathbf{U}_{1}\boldsymbol{\otimes }\mathbf{U}_{\ell
-1},\quad \mathbf{V}_{\ell }=\mathbf{V}_{1}\boldsymbol{\otimes }\mathbf{V}%
_{\ell -1},\vspace{0.05in} \\ 
\mathbf{\Sigma }_{A\ell }=\mathbf{\Sigma }_{E1}\boldsymbol{\otimes }\mathbf{%
\Sigma }_{L,\ell -1},\quad \mathbf{\Sigma }_{B\ell }=\mathbf{\Sigma }_{L1}%
\boldsymbol{\otimes }\mathbf{\Sigma }_{E,\ell -1},\vspace{0.05in} \\ 
\mathbf{\Sigma }_{L\ell }=\mathbf{\Sigma }_{A\ell }+\mathbf{\Sigma }_{B\ell
},\quad \mathbf{\Sigma }_{E\ell }=\mathbf{D}_{E\ell }%
\end{array}%
\right\} \quad \ell =2,3,\ldots \ ,  \label{SVDn1}
\end{equation}%
where, as in (\ref{U9}), the signs of the columns of $\mathbf{U}_{\ell }$
must be adjusted to coincide with their singular value being made positive
when numerical values are specified for the parameters of $\mathbf{L}_{\ell
}\left( c_{i},v_{i},y_{i}\right) $. The formulas (\ref{SVDn1}) can be
verified in the same manner as (\ref{Jordl1}). From (\ref{SVDn1}), starting
with (\ref{SVD3}) and (\ref{Sig9}), the singular values of $\mathbf{L}_{\ell
}\left( c_{i},v_{i},y_{i}\right) $ are found to be%
\begin{gather}
\mu _{n},\ 3^{\ell -1}\left\vert \phi \left( v_{1},y_{1}\right) \right\vert
\!,\ 3^{\ell -1}\left\vert \phi \left( v_{2},y_{2}\right) \right\vert
\!,\ldots ,\ 3^{\ell -1}\left\vert \phi \left( v_{\ell },y_{\ell }\right)
\right\vert \!,  \notag \\
3^{\ell -1}\left\vert \psi \left( v_{1},y_{1}\right) \right\vert \!,\
3^{\ell -1}\left\vert \psi \left( v_{2},y_{2}\right) \right\vert \!,\ldots
,\ 3^{\ell -1}\left\vert \psi \left( v_{\ell },y_{\ell }\right) \right\vert
\!,\quad \ell =1,2,\ldots \ .  \label{SVL3n}
\end{gather}%
with $n=3\ell $ and $\phi \left( v_{i},y_{i}\right) $ and $\psi \left(
v_{i},y_{i}\right) $ from (\ref{PhiPsi}). These singular values also are
valid for compound Frierson matrices and they agree with those given by Loly
and Cameron \cite{LOLY3} since the signs of $v_{i}$ and $y_{i}$ do not
affect $\left\vert \phi \left( v_{i},y_{i}\right) \right\vert $ and $%
\left\vert \psi \left( v_{i},y_{i}\right) \right\vert $. We note that $%
\mathbf{L}_{\ell }$ and $\mathbf{F}_{\ell }$ are rank $\ell +1$
matrices.\smallskip

\noindent \textbf{Natural Squares}. The FNC (\ref{Frob3}) applied to the
singular values (\ref{SVL3n}) of $\mathbf{L}_{\ell }\left(
c_{i},v_{i},y_{i}\right) $ with (\ref{PhiPsi}) leads to%
\begin{equation}
\tsum\limits_{i=1}^{\ell }\left( v_{i}^{2}+y_{i}^{2}\right) =\frac{1}{8}%
\left( 9^{2\ell }-1\right) =\tsum\limits_{i=0}^{2\ell -1}9^{i},\quad \ell
=1,2,\ldots \ .  \label{Frobn3}
\end{equation}%
Therefore, when $v_{i}$ and $y_{i}$ $\left( i=1,2,\ldots ,\ell \right) $
take the distinct values $\pm 3^{k}$ $\left( k=0,1,\ldots ,2\ell -1\right) $
in any order, (\ref{Frob3}) is satisfied and $\mathbf{L}_{\ell }$ may be
natural. Note that (\ref{Frobn3}) yields (\ref{FrobLF}) for $\ell =1$ and (%
\ref{vyst}) for $\ell =2.$ It remains to determine $c_{i}$ for possibly
natural $\mathbf{L}_{\ell }\left( c_{i},v_{i},y_{i}\right) .$ From (\ref%
{Ll+1}), (\ref{MnMagic}), and (\ref{Mun}), it follows that%
\begin{equation}
3^{-\ell }\limfunc{tr}\mathbf{L}_{\ell }=\tsum\limits_{i=1}^{\ell
}c_{i}=3^{-\ell }\mu _{\ell }=\frac{1}{2}\left( 3^{2\ell }-1\right)
=\tsum\limits_{i=0}^{2\ell -1}3^{i}
\end{equation}%
which is satisfied by taking%
\begin{equation}
c_{i}=\left\vert v_{i}\right\vert +\left\vert y_{i}\right\vert  \label{ck}
\end{equation}%
for the FNC values for $v_{i}$ and $y_{i},$ as in (\ref{ccvy}) and (\ref%
{cdmu}).

As for $\mathbf{L}_{9},\ \mathbf{A}_{9},$ and $\mathbf{B}_{9}$ of the
preceding section, $\mathbf{A}_{\ell }$ and $\mathbf{B}_{\ell }$ of (\ref%
{AB3n}) are mutually orthogonal in form and the elements of $\mathbf{L}%
_{\ell }\left( c_{i},v_{i},y_{i}\right) $ are linear combinations of $%
c_{i},v_{i},$ and $y_{i}.$ Thus, a $\mathbf{L}_{\ell }$ that meets the FNC
would contain duplicate elements if%
\begin{equation}
\tsum\limits_{i=0}^{2\ell -1}3^{i}\alpha _{i}=\tsum\limits_{i=0}^{2\ell
-1}3^{i}\beta _{i},\quad \alpha _{i},\beta _{i}=-1,0,\text{ or }1,
\end{equation}%
i.e.%
\begin{equation}
\tsum\limits_{i=0}^{2\ell -1}3^{i}\eta _{i}=0,\quad \eta _{i}=\alpha
_{i}-\beta _{i}=-2,-1,0,1,\text{ or }2.
\end{equation}%
For the critical cases of $\eta _{i}$ we find that%
\begin{equation}
\tsum\limits_{i=0}^{2\ell -1}\eta
_{i}3^{i}=3^{k}-2\tsum\limits_{i=0}^{k-1}3^{i}=1,\quad k=1,2,\ldots ,2\ell
-1,\quad \ell =1,2,\ldots \ .
\end{equation}%
Therefore, for $v_{i}$ and $y_{i}$ $\left( i=1,2,\ldots ,\ell \right) $
taking the distinct values $\pm 3^{k}$ $\left( k=0,1,\ldots ,2\ell -1\right) 
$ in any order and with $c_{i}$ from (\ref{ck}), $\mathbf{L}_{\ell }\left(
c_{i},v_{i},y_{i}\right) $ is a natural magic square and similarly for $%
\mathbf{F}_{\ell }\left( v_{i},y_{i}\right) $ when these same $v_{i}$ and $%
y_{i}$ are positive.

Since the $2\ell $ distinct values of $v_{i}$ and $y_{i}$ can be arranged in
any order and there are $2^{2\ell }$ possible assignments of $\pm $ signs
for each ordering, the number $N_{L\ell }$ of fundamental (phases excluded)
Lucas magic squares of order $3^{\ell }$ and the number $N_{F\ell }$ of
fundamental Frierson magic squares of order $3^{\ell }$ are given by%
\footnote{%
Our formula for $N_{F\ell }$ differs from that of Loly and Cameron \cite%
{LOLY3} which gives $N_{F2}=6$ (like Frierson \cite{FIER}) and $N_{F3}=90$
instead of our $360.$%
\par
{}}%
\begin{equation}
N_{L\ell }=2^{2\ell }\left( 2\ell \right) !/8,\quad N_{F\ell }=\left( 2\ell
\right) !/2.
\end{equation}%
The number $N_{SV\ell }$ of sets of singular values that can be formed for a
compound Frierson square $\mathbf{F}_{\ell }\left( v_{i},y_{i}\right) $ is
given by Loly and Cameron \cite{LOLY3} as%
\begin{equation}
N_{SV\ell }=\left( 2\ell -1\right) !!\,=1\times 3\times 5\times \cdots
\times \left( 2\ell -1\right) \!,  \label{Nsv}
\end{equation}%
based on (\ref{SVL3n}) and the form of $\phi \left( v_{i},y_{i}\right) $ and 
$\psi \left( v_{i},y_{i}\right) $ in (\ref{PhiPsi}) using a known
combinatorial result \cite{OEIS}. In view of (\ref{SVL3n}), this formula
also applies to compound Lucas magic squares.

Pertinent numerical constants for natural $\mathbf{L}_{\ell }$ and $\mathbf{F%
}_{\ell }$ are given in Table 2.

\begin{gather*}
\begin{tabular}{||c|c|c|c|c|c|c||}
\hline\hline
$%
\begin{array}{c}
\  \\ 
\ \vspace{0.01in}%
\end{array}%
\hspace{-0.2in}\ell $ & $n$ & $\mu _{\ell }$ & $N_{L\ell }$ & $N_{F\ell }$ & 
Rank & $N_{SV\ell }$ \\ \hline
$1$ & $3$ & $12$ & $1$ & $1$ & $3$ & $1$ \\ 
$2$ & $9$ & $360$ & $48$ & $12$ & $5$ & $3$ \\ 
$3$ & $27$ & $9,828$ & $5,760$ & $360$ & $7$ & $15$ \\ 
$4$ & $81$ & $265,680$ & $1,290,240$ & $20,160$ & $9$ & $105$ \\ 
$5$ & $243$ & $7,174,332$ & $464,486,400$ & $1,814,400$ & $11$ & $945$ \\ 
$6$ & $729$ & $193,709,880$ & $245,248,819,200$ & $239,500,800$ & $13$ & $%
10,395$ \\ \hline\hline
\end{tabular}
\\
\text{Table 2 - Numerical constants for Lucas and Frierson matrices.}
\end{gather*}%
\smallskip

\noindent \textbf{Commuting Matrices}. Since $\mathbf{A}_{\ell }$ and $%
\mathbf{B}_{\ell }$ of (\ref{AB3n}) commute, based on formulas developed by
Nordgren \cite{NORD3}, the following two compound Lucas magic square
matrices commute:%
\begin{eqnarray}
&&\mathbf{L}_{\ell }^{\prime }\left(
c_{1},v_{1},y_{1},c_{2},v_{2},y_{2},\ldots ,c_{\ell -1},v_{\ell -1},y_{\ell
-1},9^{\ell -1}c_{\ell },9^{\ell -1}v_{\ell },9^{\ell -1}y_{\ell }\right) \!,
\notag \\
&&\mathbf{L}_{\ell }^{\prime \prime }\left(
9c_{1},9v_{1},9y_{1},9c_{2},9v_{2},9y_{2},\ldots ,9c_{\ell -1},9v_{\ell
-1},9y_{\ell -1},c_{\ell },v_{\ell },y_{\ell }\right) \!,  \label{LAB}
\end{eqnarray}%
where%
\begin{gather}
v_{1},y_{1},v_{2},y_{2},\ldots ,v_{\ell -1},y_{\ell -1}\Rightarrow \pm 1,\pm
3,\pm 9,\ldots ,\pm 3^{2\ell -3}\text{ in any order,}  \notag \\
v_{\ell },y_{\ell }\Rightarrow \pm 1,\pm 3\text{ in any order.}
\end{gather}%
Similar equations form commuting Frierson matrices (with positive
parameters), e.g.%
\begin{equation}
\left[ \mathbf{F}_{2}\left( 1,3,27,9\right) ,\mathbf{F}_{2}\left(
9,27,3,1\right) \right] =\left[ \mathbf{F}_{9G}\mathbf{R}_{9},\mathbf{F}_{9J}%
\right] =\mathbf{0,}  \label{CommEx}
\end{equation}%
as previously determined by (\ref{F9comm}). Other pairs of commuting $%
\mathbf{L}_{\ell }$'s and $\mathbf{F}_{\ell }$'s can be formed from $\mathbf{%
A}_{\ell }$ and $\mathbf{B}_{\ell }$ in a similar manner. In addition, two
different commuting $\mathbf{L}_{\ell }$'s or $\mathbf{F}_{\ell }$'s can be
used to form two different $\mathbf{A}_{\ell }$'s and $\mathbf{B}_{\ell }$'s
for $\mathbf{L}_{\ell +1}^{\prime }$ and $\mathbf{L}_{\ell +1}^{\prime
\prime }$ which commute, provided that the $v_{i}$ and $y_{i}$ come from $%
\pm 1,\pm 3,\pm 9,\ldots ,\pm 3^{2\ell +1}$ so that the resulting $\mathbf{L}%
_{\ell +1}^{\prime }$ and $\mathbf{L}_{\ell +1}^{\prime \prime }$ are
natural magic squares. A similar construction applies to Frierson squares,
e.g.%
\begin{equation}
\left[ \mathbf{F}_{3}\left( 1,3,27,9,81,243\right) ,\mathbf{F}_{3}\left(
9,27,3,1,81,243\right) \right] =\mathbf{0}
\end{equation}%
which can easily be verified.

\section{Conclusion}

Lucas' parameterization of order-3 magic square matrices \cite{LUCA} can be
compounded sequentially by a known procedure to produce parameterized Lucas
magic square matrices of order $3^{\ell }\ \left( \ell =2,3,\ldots \right) .$
Formulas from \cite{NORD3} determine the matrices in their Jordan canonical
form and singular value decomposition. The Frobenius norm provides a
necessary condition for a magic square matrix to be natural. Values for the
parameters that produce natural compound Lucas magic square matrices follow
from formulas for their singular values and the Frobenius norm condition
together with verification that they have no duplicate elements. Frierson's
parameterized magic squares \cite{FIER,LOLY3} are a special case of compound
Lucas squares and our results apply to them as well.

We determine the number of fundamental (phases excluded) Lucas and Frierson
natural magic squares that can be constructed by the compounding procedure.
We find 48 specific fundamental Lucas squares of order 9. We find six new
specific fundamental Frierson squares of order 9 in addition to the six
given by Frierson. Also, we find commuting pairs of compound Lucas matrices
and formulas for matrix powers of order-3 and order-9 Lucas matrices.%
\footnote{%
The derivations in this article were carried out using Maple$^{\copyright }$
and MuPAD$^{\copyright }$ in Scientific WorkPlace$^{\copyright }\!$.}

\end{document}